\documentclass[12pt]{report}
\usepackage{latexsym}
\newcommand{\hs}{\hspace{0.3cm}}
\newcommand{\insp}{\hspace*{1cm}}

\newcommand{\ds}{\displaystyle}
\newcommand{\ts}{\textstyle}

\usepackage{epsfig}

\newcommand{\ind}{\mbox{ind}}
\newcommand{\Aut}{\mbox{Aut}}
\newcommand{\ad}{\mbox{ad}}
\newcommand{\Gr}{\mbox{Gr}}

\newcommand{\reg}{\mbox{\scriptsize{reg}}}

\newcommand{\rank}{\mbox{rank}}

\newcommand{\Der}{\mbox{Der}}

\newcommand{\IC}{I\!\!\!\!C}

\newcommand{{\CB}}{\cal B}

\newcommand{{\CC}}{\cal C}
\def\subs{\mathop{\subset}}

\font\myfont=cmfi10 scaled \magstep 1

\def\picture#1 by #2 (#3){
\vbox to #2{
\hrule width #1 height 0pt depth 0pt
\vfill
\special{picture #3}}}

\begin{document}
\begin{center}
{\bf \Large On Commutative Polarizations}\\
\ \\
{\bf \large Alexander G. Elashvili\footnote{Partially supported by the CRDF
grant RM1-2088}}\\
{\it Razmadze Mathematics Institute, M. Aleksidze Str.1, 380093 Tbilisi,\\
Republic of Georgia\\
E-mail address: alela@rmi.acnet.ge}\\
\ \\
and\\
\ \\
{\bf \large Alfons I. Ooms}\\
{\it University of Limburg, Mathematics Department, 3590 Diepenbeek, Belgium\
E-mail address: alfons.ooms@luc.ac.be}
\end{center}

\

\pagestyle{plain}
\begin{center}
{\bf INTRODUCTION}
\end{center}
Let $L$ be a finite-dimensional Lie algebra over a field $k$ of characteristic
zero and let $U(L)$ be its enveloping algebra with quotient division ring
$D(L)$.  Let $P$ be a commutative Lie subalgebra of $L$.  In [O2] the necessary
and sufficient condition on $P$ was given in order for $D(P)$ to be a maximal
(commutative) subfield of $D(L)$.  In particular, this condition is satisfied if
$P$ is a commutative polarization (CP) with respect to any regular $f \in
L^\ast$ and the converse holds if $L$ is ad-algebraic.  The purpose of this
paper is to study Lie algebras admitting these CP's and to demonstrate their
widespread occurrence.\\
First we have the following characterisation if $L$ is completely solvable: $P$
is a CP of $L$ if and only if there exists a descending chain of Lie subalgebras
\begin{eqnarray*}
L = L_n \supset \ldots \supset L_{j+1}\supset L_j \supset \ldots \supset L_p = P
\end{eqnarray*}
such that \hs $\dim L_j = j$ \hs with increasing index, i.e. \hs
$i(L_j) = i(L_{j+1}) + 1$, \hs $j : p,\ldots, n-1$ (Theorem 1.11).
In low dimension this phenomenon appears frequently. In fact, in a
case by case study of indecomposable nilpotent Lie algebras of
dimension at most seven we discover that Lie algebras without CP's
are rather exceptional: 1 (out of 9) in dimension at most 5; 3
(out of 22) in dimension 6 and 26 (out of 130) in dimension 7.
These will be listed in section 3, in which we also prove that
nonabelian Lie algebras having a nondegenerate, invariant bilinear
form do not admit any CP (Theorem 3.2).

\

Suppose $k$ is algebraically closed.  Then for a Lie algebra $L$
to admit a CP $P$ has the following advantage: in $U(L)$ the
primitive ideals $I(f)$, with regular $f \in L^\ast$, can all be
constructed using the same polarization $P$, since $I(f)$ is the
kernel of the (twisted) induced representation \hs $\sigma =
\ind^\sim (f|_P,L)$ \hs [D, 10.3.4].  If in addition $P$ is an
ideal of $L$ (a so called CP-ideal of $L$) then the representation
$\sigma$ is irreducible (in the completely solvable case $P$ even
turns out to be a Vergne polarization). Also, the semi-center
$Sz(U(L))$ of $U(L)$ is contained in $U(P)$ (Corollary 4.4).
Moreover, a standard technique using Grassmannians shows that if
$L$ is solvable with a CP, then it also has a CP-ideal (Theorem
4.1).

\

In section 5, we look for CP-ideals in some Frobenius Lie algebras
(i.e. Lie algebras of index zero [O1]).  For instance, let $x \in
L$ be a principal nilpotent element of a semi-simple Lie algebra
$L$ with centralizer $P$.  Then the normalizer $F$ of $P$ is a
Frobenius Lie algebra by a recent result of Panyushev [P2], in
which $P$ is a CP-ideal (Theorem 5.7).  Next, let $A$ be a finite
dimensional associative algebra over $k$ with a unit.  A becomes a
Lie algebra {\myfont g} for the Lie bracket \hs $[a,b] = ab-ba$,
\hs $a,b \in A$ \hs and $V = A$ becomes a {\myfont g}-module by
left multiplication.  Consider the semi-direct product \hs $L = \
${\myfont g}$\ \oplus\ V$. \hs Then the following are equivalent
(Proposition 5.6):
\begin{itemize}
\item[(1)] $A$ is a Frobenius algebra
\item[(2)] $L$ is a Frobenius Lie algebra
\item[(3)] $V$ is a CP-ideal of $L$
\item[(4)] $D(V)$ is a maximal subfield of $D(L)$.
\end{itemize}
A similar result can be obtained if $A$ is a finite dimensional
left symmetric algebra (Example 5.4) or if $A$ is a finite
dimensional simple Novikov algebra over $k$, \hs char$(k) = p >
2$.

\

CP-ideals also occur naturally in the nilradical $N$ of any
parabolic Lie subalgebra of a simple Lie algebra $L$ of type $A_r$
or $C_r$.  As a bonus we
obtain an explicit formula for the index $i(N)$ of $N$ (Theorem 6.2).\\
Finally, section 7 deals with some CP-preserving extensions.
\newpage
\begin{center}
{\bf 1. PRELIMINARIES AND GENERAL RESULTS}
\end{center}

Let $L$ be a Lie algebra over a field $k$ of characteristic zero with basis \hs
$x_1,\ldots,x_n$. \hs Let \hs $f \in L^\ast$ \hs and consider the alternating
bilinear form $B_f$ on $L$ sending \hs $(x,y)$\hs into \hs $f([x,y])$.\hs For
any subset $A$ of $L$ we denote by $A^\perp$ or $A^f$ the subspace
\begin{eqnarray*}
\{x \in L \mid f([x,a]) = 0\hs \mbox{for all} \hs a \in A \}
\end{eqnarray*}
We also put \hs $L(f) = L^\perp$ \hs and \hs $i(L) = \min\limits_{f \in L^\ast}
\dim L(f)$, the index of $L$.  Note that $L(f)$ is a Lie subalgebra of $L$
containing the center $Z(L)$ of $L$. We recall from [D, 1.14.13] that
\begin{eqnarray*}
i(L) = \dim L - \rank_{R(L)} ([x_i, x_j])
\end{eqnarray*}
where $R(L)$ is the quotient field of the symmetric algebra $S(L)$ of $L$.  In
particular, \hs $\dim L - i(L)$ \hs is an even number.\\
Furthermore, $f$ is called regular if \hs $\dim L(f) = i(L)$. \hs
It is well- known that the set $L_{\reg}^\ast$ of all regular
elements of $L^\ast$ is an open dense subset of $L^\ast$ for the
Zariski topology.

\

{\bf DEFINITION 1.1} \hs [D, 1.12.7] A Lie subalgebra $P$ of $L$
is called a polarization w.r.t. \hs $f \in L^\ast$ \hs if \hs
$f([P,P]) = 0$ \hs and \hs $\dim P = \frac{1}{2} (\dim L + \dim
L(f))$, \hs in other words $P$ is a maximal totally isotropic
subspace of $L$ (equipped with $B_f$).  If in addition $P$ is
commutative then $f$ is regular by the following observation.

\

{\bf LEMMA 1.2} \hs (see Theorem 14 of [O2]). Let $P$ be a
commutative Lie subalgebra of $L$; \hs $h_1,\dots, h_m$\hs a basis
of $P$ and \hs $x_1,\ldots, x_n$ \hs a basis of $L$.  Then the
following conditions are equivalent:
\begin{itemize}
\item[(a)] $\dim P = \frac{1}{2} (\dim L + i(L))$, \hs i.e. $P$ is a CP
(commutative polarization) of $L$ w.r.t. each $f \in L_{\reg}^\ast$.
\item[(b)] $P = P^f$ \hs w.r.t. some \hs $f \in L^\ast$ \hs (such an $f$ is
necessarily regular)
\item[(c)] $\rank_{R(L)} ([h_i, x_j]) = \dim L - \dim P$
\end{itemize}

{\bf LEMMA 1.3} \hs Let $P$ and $M$ be Lie subalgebras of $L$ such that \hs $P
\subset M \subset L$. \hs Then the following conditions are equivalent:
\begin{itemize}
\item[(1)] $P$ is a CP of $L$
\item[(2)] $P$ is a CP of $M$ and \hs $i(M) = i(L) + \dim L - \dim M$.
\end{itemize}
Under these conditions the following hold:
\begin{eqnarray*}
f \in L_{\reg}^\ast \hs \Rightarrow \hs f|_M \in M_{\reg}^\ast
\end{eqnarray*}

{\it Proof.} \hs (1) $\Rightarrow$ (2). Take any \hs $f \in L_{\reg}^\ast$.\hs
Then \hs $P = P^f$.  Put \hs $g = f|_M \in M^\ast$.  \hs W.r.t. $B_g$ we have:
\begin{eqnarray*}
P^g &=& \{x \in M\mid g([x,P]) = 0\} = \{x \in M \mid f([x,P]) = 0\}\\
&=& M \cap P^f = M\cap P = P.
\end{eqnarray*}
Hence $P$ is a CP of $M$ and \hs $g \in M_{\reg}^\ast$\hs by Lemma 1.2.  In
particular,
\begin{eqnarray*}
\ts\frac{1}{2} (\dim M + i(M)) = \dim P = \frac{1}{2} (\dim L + i(L))
\end{eqnarray*}
Consequently, \hs $i(M) = i(L) + \dim L - \dim M$.\\
(2) $\Rightarrow$ (1). \hs $P$ is commutative and
\begin{eqnarray*}
\dim P &=& \ts\frac{1}{2} (\dim M + i(M))\\
&=& \ts\frac{1}{2} (\dim M + i(L) + \dim L - \dim M)\\
&=& \ts\frac{1}{2} (\dim L + i(L))
\end{eqnarray*}
Hence, $P$ is a CP of $L$.

\

The following is a direct application of [D, Lemma 1.12.2].

\

{\bf LEMMA 1.4.} \hs Let $M$ be a Lie subalgebra of $L$ of
codimension one.  Let \hs $f \in L^\ast$\hs and put \hs $g = f|_M
\in M^\ast$.  \hs Then we distinguish two cases:
\begin{itemize}
\item[(i)] If \hs $L(f) \subset M$ \hs then $L(f)$ is a hyperplane
in $M(g)$. \item[(ii)] If \hs $L(f) \not\subset M$ \hs then \hs
$M(g) = L(f) \cap M$\hs is a hyperplane in $L(f)$.
\end{itemize}

\

{\bf REMARK 1.5.} \hs In [O2] we introduced the notion of the Frobenius
semiradical $F(L)$ of a Lie algebra $L$, namely
\begin{eqnarray*}
F(L) = \sum\limits_{f \in L_{\reg}^\ast} L(f)
\end{eqnarray*}
This is a characteristic ideal of $L$ containing the center $Z(L)$
of $L$.  It seems to play a natural role in the study of
commutative polarizations.  For instance if $L$  admits a CP $P$,
then \hs $F(L) \subset P$ \hs and hence is commutative [O2,
p.710].

\

{\bf PROPOSITION 1.6.} \hs Let $M$ be a Lie subalgebra of $L$ of
codimension one, \hs $f \in L^\ast$\hs and \hs $g = f|_M \in
M^\ast$. \hs Then we have:
\begin{itemize}
\item[(1)] either \hs $i(M) = i(L) + 1$ \hs or \hs $i(M) = i(L) -1$
\item[(2)] $\left\{\begin{array}{l} f \in L_{\reg}^\ast\\i(M) = i(L) +
1\end{array}\right. \hs \Leftrightarrow \hs \left\{\begin{array}{l} g \in
M_{\reg}^\ast\\L(f) \subset M\end{array}\right.$
\item[(3)] $\left\{\begin{array}{l} f \in L_{\reg}^\ast\\L(f) \not\subset
M\end{array}\right.\hs \Leftrightarrow \insp \hs \left\{\begin{array}{l} g\in
M_{\reg}^\ast\\i(M) = i(L) -1\end{array}\right.$
\item[(4)] $i(M) = i(L) + 1 \hs \Leftrightarrow \hs F(L) \subset M$
\item[(5)] Suppose \hs $i(M) = i(L) + 1$ \hs and let $P$ be a Lie subalgebra of
$M$.  Then
\begin{eqnarray*}
P\hs \mbox{is a CP of}\ L \hs \Leftrightarrow\hs P\ \mbox{is a CP of}\ M
\end{eqnarray*}
\item[(6)] Suppose \hs $i(M) = i(L) -1$. \hs If $H$ is a CP (respectively a CP-
ideal) of $L$, then \hs $H \cap M$\hs is a CP (resp. a CP-ideal) of $M$ and \hs
$\dim(H\cap M) = \dim H-1$.
\end{itemize}

{\it Proof.}
\begin{itemize}
\item[(1)] Choose \hs $\varphi \in L_{\reg}^\ast$ \hs such that \hs $\gamma =
\varphi|_M \in M_{\reg}^\ast$. \hs Suppose \hs $L(\varphi) \subset M$ \hs then
\begin{eqnarray*}
i(M) = \dim M(\gamma) = \dim L(\varphi) + 1 = i(L) + 1
\end{eqnarray*}
by (i) of Lemma 1.4.  On the other hand, if \hs $L(\varphi) \not\subset M$ \hs
then
\begin{eqnarray*}
i(M) = \dim M(\gamma) = \dim L(\varphi) -1 = i(L) - 1
\end{eqnarray*}
by (ii) of Lemma 1.4.
\item[(2)] $\Rightarrow :$\\
Suppose $L(f) \not\subset M$.  By (ii) of Lemma 1.4
\begin{eqnarray*}
i(M) \leq \dim M(g) = \dim L(f) - 1 = i(L) - 1
\end{eqnarray*}
Contradiction.  Therefore \hs $L(f) \subset M$. \hs Hence,
\begin{eqnarray*}
i(M) - 1 = i(L) = \dim L(f) = \dim M(g) - 1
\end{eqnarray*}
by (i) of Lemma 1.4. Hence \hs $i(M) = \dim M(g)$, \hs i.e. \hs $g \in
M_{\reg}^\ast$.\\
$\Leftarrow :$\\
By (i) of Lemma 1.4 \hs $L(f) \subset M$ \hs implies that
\begin{eqnarray*}
i(L) \leq \dim L(f) = \dim M(g) - 1 = i(M) -1
\end{eqnarray*}
So, \hs $i(M) \geq i(L) + 1$. \hs By (1), \hs $i(M) = i(L) + 1$ \hs and
therefore \\
$i(L) = \dim L(f)$, \hs i.e. \hs $f \in L_{\reg}^\ast$.
\item[(3)] $\Rightarrow :$\\
$L(f) \not\subset M$ \hs implies that
\begin{eqnarray*}
i(M) \leq \dim M(g) = \dim L(f) - 1 = i(L) - 1
\end{eqnarray*}
by (ii) of Lemma 1.4.  Hence, by (1), \hs $i(M) = i(L) -1$ \hs which forces \\
$i(M) = \dim M(g)$,\hs i.e. \hs $g \in M_{\reg}^\ast$.\\
$\Leftarrow :$\\
Since \hs $i(M) \neq i(L) + 1$ \hs it follows from (2) that \hs $L(f)
\not\subset M$. \hs Hence,
\begin{eqnarray*}
i(L) - 1 = i(M) = \dim M(g) = \dim L(f) - 1
\end{eqnarray*}
by (ii) of Lemma 1.4.  Consequently, \hs $\dim L(f) = i(L)$, \hs i.e. \hs $f \in
L_{\reg}^\ast$.
\item[(4)] $\Rightarrow$ follows from (2).\\
$\Leftarrow$ Choose \hs $f \in L_{\reg}^\ast$ \hs such that \hs $g = f|_M \in
M_{\reg}^\ast$.  \hs Then \\
$L(f) \subset F(L) \subset M$. \hs Using (2) it follows that \hs $i(M) = i(L)
+1$.
\item[(5)] Clearly, \hs $i(M) = i(L) + \dim L - \dim M$. Now use Lemma 1.3.
\item[(6)] Suppose \hs $i(M) = i(L) - 1$.\hs Hence, by Lemma 1.3 \hs $H
\not\subset M$. \hs Then \\
$\dim (H \cap M) = \dim H -1$. \hs $H \cap M$ \hs is abelian and
\begin{eqnarray*}
\dim (H \cap M) = \ts\frac{1}{2} (\dim L + i(L)) - 1 = \ts\frac{1}{2} (\dim M +
i(M))
\end{eqnarray*}
Consequently, \hs $H \cap M$\hs is a CP (resp. a CP-ideal) of $M$.
\end{itemize}

\newpage
{\bf EXAMPLES 1.7.}
\begin{itemize}
\item[(1)] Let $E$ be a nonzero endomorphism of an $n$-dimensional vector space
$V$ over $k$.  Consider the Lie algebra \hs $L = kE \oplus V$ \hs with Lie
brackets \hs $[E,v] = Ev$ \hs and in which $V$ is a commutative ideal.  $L$ is
solvable and \hs $i(L) = n-1$. \hs Clearly, \hs $i(V) = n = i(L) + 1$ \hs and
$V$ is a CP-ideal of $L$ by (5) of Proposition 1.6.
\item[(2)] Let $L$ be a Frobenius Lie algebra (i.e. \hs $i(L) = 0$) and $M$ a
Lie subalgebra of $L$ of codimension one. Then \hs $i(M) = 1\hs (=i(L) +1)$.
\item[(3)] Let $M$ be a Lie subalgebra of codimension one in a nonabelian Lie
algebra $L$ with \hs $F(L) = L$. \hs Then, \hs $i(M) = i(L) - 1$ \hs and $L$
does not have any CP's (by Proposition 1.6 and Remark 1.5).  For instance, let
$L$ be the diamond Lie algebra with basis \hs $t, x, y, z$ \hs and nonvanishing
brackets \hs $[t,x] = -x \hs [t,y] = y$ \hs and \hs $[x,y] = z$.\hs Clearly, \hs
$i(L) = 2$ \hs and \hs $M = [L,L] = \langle x,y,z\rangle$ \hs is an ideal of
codimension one in $L$ with \hs $i(M) = 1$. \hs Put \hs $f = x^\ast \in
L_{\reg}^\ast$\hs and \hs $g = f|_M \in M^\ast$.\hs Then, \hs $L(f) = \langle
y,z\rangle \subset M$, \hs $i(M) = i(L) -1$\hs and \hs $g \notin M_{\reg}^\ast$.
\hs Also, \hs $P_1 = \langle y,z\rangle$\hs is a CP of $M$.  But there is no CP
\hs $P$ of $L$ such that \hs $P \cap M = P_1$ \hs (in fact $L$ does not admit
any CP since \hs $F(L) = L$).\\
See also Theorem 3.2 and (2) of Examples 3.3.
\end{itemize}

\

{\bf DEFINITION 1.8} \hs A Lie algebra $L$ is called square integrable if \hs
$L(f) = Z(L)$ \hs for some\hs $f \in L^\ast$, \hs i.e. \hs $i(L) = \dim Z(L)$.\\
In the nilpotent case these Lie algebras are precisely the Lie
algebras of simply connected Lie groups admitting square
integrable representations [MW, p.450-453].

\

{\bf PROPOSITION 1.9} \hs Let $L$ be a Lie algebra having an element $u \in L$
such that its centralizer \hs $M = C(u)$ \hs has codimension one in $L$.  Then
we have
\begin{itemize}
\item[(i)] $i(M) = i(L) + 1$ \item[(ii)] $L$ has a CP if and only
if $M$ has a CP \item[(iii)] If $L$ is square integrable then so
is $M$.
\end{itemize}

\

{\bf REMARK 1.10.} \hs Note that $C(u)$ is an ideal of codimension
one of $L$ if either $u$ is a noncentral semi-invariant of $L$
(i.e. for a suitable \hs $\lambda \in L^\ast \backslash \{0\} :
\hs [x,u] = \lambda (x) u,\hs x \in L$) or \hs $[u,L]$ \hs is a
one dimensional subspace of the center $Z(L)$ (such an $u$ always
exists if $L$ is nilpotent and \hs $\dim Z(L) = 1 < \dim L$).  In
that situation, if $L$ has a CP-ideal then the same holds for
$C(u)$.

\

\

{\it Proof of the proposition.}
\begin{itemize}
\item[(i)] Take \hs $x \in L \backslash C(u)$ \hs and choose \hs $f \in L^\ast$
\hs such that $f|_M$ is regular and such that \hs $f([x,u]) \neq 0$. \hs Then
\hs $C(u) = u^f$ \hs (since both have the same dimension and \hs $C(u) \subset
u^f$).  Then \hs $L(f) = L^f \subset u^f = M$. \hs It follows by (2) of
Proposition 1.6 that \hs $i(M) = i(L) + 1$ \hs and \hs $f \in L_{\reg}^\ast$.
\item[(ii)] First, let $P$ be a commutative Lie subalgebra of $M$.  Then,
\begin{center}
$P$ is a CP of $L$ if and only if $P$ is a CP of $M$
\end{center}
by (5) of Proposition 1.6.  Next, let $P$ be a CP of $L$ such that \hs $P
\not\subset M$.\hs Then \hs $\dim (P \cap M) = \dim P -1$ \hs and \hs $u \notin
P\cap M$\hs (otherwise \hs $[u,P] = 0$ \hs and thus \hs $P \subset C(u) = M$).\\
Finally, \hs $P_1 = (P \cap M) \oplus ku$ \hs is a CP of $M$ since it is
commutative and
\begin{eqnarray*}
\dim P_1 = \dim P = \ts\frac{1}{2} (\dim L + i(L)) = \ts\frac{1}{2}(\dim M +
i(M))
\end{eqnarray*}
\item[(iii)] Clearly, \hs $Z(L) \subset C(u) = M$ \hs and \hs $u \in Z(M)
\backslash Z(L)$.\hs Hence, \\
$Z(L) \oplus ku \subset Z(M)$. \hs Therefore,
\begin{eqnarray*}
i(M) \geq \dim Z(M) \geq \dim Z(L) + 1 = i(L) + 1
\end{eqnarray*}
As \hs $i(M) = i(L) + 1$ \hs we may conclude that \hs $i(M) = \dim
Z(M)$, \hs i.e. $M$ is square integrable.
\end{itemize}

\

{\bf THEOREM 1.11} \hs Let $P$ be a commutative Lie subalgebra of a completely
sol\-vable Lie algebra $L$.  Then the following conditions are equivalent:
\begin{itemize}
\item[(1)] $P$ is a CP (resp. CP-ideal) of $L$.
\item[(2)] There exists a descending series of Lie subalgebras (resp. ideals) of
$L$.
\end{itemize}
\begin{eqnarray*}
L = L_n \supset \ldots \supset L_{j+1} \supset L_j \supset \ldots \supset L_p =
P
\end{eqnarray*}
$\dim L_j = j$, \hs with increasing index (i.e. \hs $i(L_j) =
i(L_{j+1}) + 1$).

\

{\it Proof.} \hs Let $P$ be a Lie subalgebra (resp. ideal) of $L$.
\hs $P$ (resp. $L$) acts on the quotient space \hs $L/P$.\hs
Application of Lie's theorem to this action shows the exis\-tence
of Lie subalgebras (resp. ideals) $L_j$ of $L$ such that \hs $L =
L_n \supset \ldots \supset L_p = P$ \hs with \hs
$\dim L_j = j$.\\
(1) $\Rightarrow$ (2).\hs Now suppose $P$ is a CP of $L$.\\
Then, by Lemma 1.3, $P$ is also a CP for each $L_j$ and
\begin{eqnarray*}
i(L_j) &=& i(L) + (n-j) = i(L) + (n-(j+1)) + 1\\
&=& i(L_{j+1}) + 1
\end{eqnarray*}
(2) $\Rightarrow$ (1)\\
By induction on $j$ we show that $P$ is a CP of $L_j$.  This is
trivial for \hs $j = p$. \hs Next, let \hs $j \geq p+1$. \hs Then
$P$ is a CP of $L_{j-1}$ and also of $L_j$ since \hs $i(L_{j-1}) =
i(L_j) +1$ \hs by (5) of Proposition 1.6.

\

{\bf COROLLARY 1.12.} \hs Let $L$ be a completely solvable
Frobenius Lie algebra of dimension $2n$ having a CP $P$.  Then $L$
can be obtained from the $n$- dimensional abelian Lie algebra $P$
with $n$ successive extensions as described in Theorem 1.10.

\

{\bf LEMMA 1.13.} \hs Let $P$ be a CP (resp. a CP-ideal) of a Lie
algebra $L$, $A$ an ideal of $L$ contained in $P$ and \hs $f \in
L_{\reg}^\ast$ \hs such that \hs $f(A) = 0$. \hs Then \hs $P/A$
\hs is a CP (resp. a CP-ideal) of the Lie algebra \hs $L/A$ \hs
and
\begin{eqnarray*}
i(L/A) = i(L) - \dim A
\end{eqnarray*}

{\it Proof.} \hs Let \hs $\varphi : L \rightarrow L/A$ \hs be the quotient
homomorphism. As \hs $f(A) = 0$\hs there is a \hs $g \in (L/A)^\ast$ \hs such
that \hs $g \circ \varphi = f$. \hs Clearly, \hs $P/A$ \hs is an abelian Lie
subalgebra (resp. ideal) of \hs $L/A$. \hs It suffices to show that \hs $(P/A)^g
= P/A$.
\begin{eqnarray*}
(P/A)^g &=& \{\varphi (x) \in L/A \mid g([\varphi(x),\varphi(P)]) = 0, x \in
L\}\\
&=& \varphi(\{x \in L \mid f([x,P]) = 0\}\\
&=& \varphi (P^f) = \varphi(P) = P/A
\end{eqnarray*}
as \hs $P^f = P$.  So, by Lemma 1.2 \hs $P/A$ \hs is a CP (resp. CP-ideal) of
\hs $L/A$ \hs and \hs $g \in (L/A)_{\reg}^\ast$. \hs Therefore, \hs $\dim P/A =
\frac{1}{2} (\dim L/A + i(L/A))$ \hs and
\begin{eqnarray*}
i(L/A) &=& 2 \dim P/A - \dim L/A\\
&=& 2(\dim P - \dim A) - (\dim L - \dim A)\\
&=& (2 \dim P - \dim L) - \dim A = i(L) - \dim A
\end{eqnarray*}

\

\

\begin{center}
{\bf 2. CP'S IN SQUARE INTEGRABLE NILPOTENT LIE ALGEBRAS}
\end{center}
The following lemma is easy to verify.

\

{\bf LEMMA 2.1.} \hs Suppose $L$ is a direct product of Lie
algebras; \hs $L = L_1 \times L_2$. \hs Then we have the
following:
\begin{itemize}
\item[(1)] $i(L) = i(L_1) + i(L_2)$ \hs and \hs $Z(L) = Z(L_1) \times Z(L_2)$.
\item[(2)] $L$ is square integrable if and only if the same holds for $L_1$ and
$L_2$.
\item[(3)] $L$ has a CP (resp. CP-ideal) if and only if the same holds for $L_1$
and $L_2$.
\end{itemize}

\

{\bf PROPOSITION 2.2.} \hs Let $L$ be a square integrable
nilpotent Lie algebra over $\IC$, of dimension $n$ at most seven.
Then $L$ admits a CP-ideal.

\

{\it Proof.} \hs By Lemma 2.1 we may assume that $L$ is
indecomposable.  In particular, \hs $1 \leq \dim Z(L) = i(L) <
\dim L$.

\

We now distinguish the following cases:
\begin{itemize}
\item[(1)] $i(L) = 1$. \hs Then $n$ is 3, 5 or 7.  Let $m$ be the maximum
dimension of all abelian ideals of $L$.  Then by [Mo, p.161] and [O2, p.706] we
have the following inequalities:
\begin{eqnarray*}
\frac{1}{2}(\sqrt{8n + 1} - 1) \leq m \leq \frac{1}{2}(\dim L + i(L)) =
\frac{1}{2}(n+1)
\end{eqnarray*}
This implies that \hs $m = \frac{1}{2} (\dim L + i(L))$ \hs in case $n = 3, 5$
or 7, showing the existence of a CP-ideal in $L$.
\item[(2)] $i(L) = 2$. \hs Then \hs $n = 6$ \hs (The case \hs $n = 4$ \hs does
not occur since $L$ is indecomposable).  We select from Morozov's classification
of 6-dimensional nilpotent Lie algebras those that are indecomposable, square
integrable and of index 2; in each \hs $\{e_1,\ldots, e_6\}$ \hs is a basis of
$L$.  The numbe\-ring is Morozov's [Mo, p.168].
\begin{itemize}
\item[4.] $[e_1, e_2] = e_5$, \hs $[e_1, e_3] = e_6$, \hs $[e_2, e_4] = e_6$
\item[5.] $[e_1, e_3] = e_5$, \hs $[e_1, e_4] = e_6$, \hs $[e_2, e_4] = e_5$ \hs
$[e_2, e_3] = \gamma e_6$, \hs $\gamma \neq 0$
\item[6.] $[e_1, e_2] = e_6$, \hs $[e_1, e_3] = e_4$, \hs $[e_1, e_4] = e_5$,
\hs $[e_2, e_3] = e_5$
\item[7.] $[e_1, e_3] = e_4$, \hs $[e_1, e_4] = e_5$, \hs $[e_2, e_3] = e_6$
\item[8.] $[e_1, e_2] = e_3 + e_5$, \hs $[e_1, e_3] = e_4$, \hs $[e_2, e_5] =
e_6$
\item[9.] $[e_1, e_2] = e_3$, \hs $[e_1, e_3] = e_4$, \hs $[e_1, e_5] = e_6$,
\hs $[e_2, e_3] = e_6$
\item[10.] $[e_1, e_2] = e_3$, \hs $[e_1, e_3] = e_5$, \hs $[e_1, e_4] = e_6$,
\hs $[e_2, e_4] = e_5$, \hs $[e_2, e_3] = \gamma e_6$, \hs $\gamma \neq 0$
\item[11.] $[e_1, e_2] = e_3$, \hs $[e_1, e_3] = e_4$, \hs $[e_1, e_4] = e_5$,
\hs $[e_2, e_3] = e_6$
\end{itemize}
In each one of these, \hs $P = \langle e_3, e_4, e_5, e_6\rangle$ \hs is a CP-
ideal, since $P$ is an abelian ideal and \hs $\dim P = 4 = \frac{1}{2} (\dim L +
i(L))$.
\item[(3)] $i(L) = 3$.  \hs Then \hs $n = 7$ (The case \hs $n = 5$ \hs does not
occur since $L$ is indecomposable).\\
We have the following possibilities according to Seeley's classification of 7-
dimensional nilpotent Lie algebras.  We maintain the same notation as in [See].
In particular \hs $\{a, b, c, d, e, f, g\}$ \hs is a basis of $L$.  In each case
we exhibit a commutative ideal $P$ of dimension 5 \hs $(= \frac{1}{2} (\dim L +
i(L)))$.\\
In the following 3 Lie algebras we take \hs $P = \langle a, d, e, f, g\rangle$\\
$3\ 7_B$: $[a,b] = e$, \hs $[b,c] = f$, \hs $[c,d] = g$\\
$3\ 7_C$: $[a,b] = e$, \hs $[b,c] = f$, \hs $[c,d] = e$, \hs $[b,d] = g$\\
$3\ 7_D$: $[a,b] = e$, \hs $[b,d] = g$, \hs $[c,d] = e$, \hs $[a,c] = f$\\
In the following 3 we take \hs $P = \langle c,d,e,f,g\rangle$\\
$3,5,7_A$: $[a,b] = c$, \hs $[a,c] = e$, \hs $[a,d] = g$, \hs $[b,d] = f$\\
$3,5,7_B$: $[a,b] = c$, \hs $[a,c] = e$, \hs $[a,d] = g$, \hs $[b,c] = f$\\
$3,5,7_C$: $[a,b] = c$, \hs $[a,c] = e$, \hs $[a,d] = g$, \hs $[b,c] = f$, \hs
$[b,d] = e$
\end{itemize}

\

{\bf REMARK 2.3} \hs Among the Lie algebras described in
Proposition 2.2 there is one which is characteristically
nilpotent, namely $1,2,4,5,7_N$ with basis $\{a,b,c,d,e,f,g\}$
and nonzero brackets: $[a,b] = c$, \hs $[a,c] = d$, \hs $[a,d] =
g$, \hs $[a,e] = f$, \hs $[a,f] = g$, \hs $[b,c] = e$, \hs $[b,d]
= f$, \hs $[b,e] = \xi g$, \hs $[b,f] = g$, \hs $[c,d] = g$, \hs
$[c,e] = -g$ \hs with \hs $\xi \neq 0, 1$. \hs [See, p.493].  In
this case take \hs $P = \langle d, e, f, g\rangle$.

\

\begin{center}
{\bf 3. LIE ALGEBRAS WITHOUT CP'S}
\end{center}
First we want to show that the restriction on the dimension in
Proposition 2.2 cannot be removed.

\

{\bf EXAMPLES 3.1}
\begin{itemize}
\item[(i)] Let $L$ be the 8-dimensional Lie algebra over $k$ with basis \hs
$\{e_1,\ldots, e_8\}$ \hs and nonvanishing brackets: \hs $[e_1, e_2] = e_5$, \hs
$[e_1,e_3] = e_6$, \hs $[e_1,e_4] = e_7$, \hs $[e_1, e_5] = -e_8$, \hs $[e_2,
e_3] = e_8$, \hs $[e_2, e_4] = e_6$, \hs $[e_2, e_6] = -e_7$, \hs $[e_3, e_4] =
-e_5$, \hs $[e_3, e_5] = -e_7$, \hs $[e_4, e_6] = -e_8$.\\
$L$ is characteristically nilpotent [DL].  $L$ is also square
integrable of index 2, but it does not admit a CP-ideal (and not
any CP's either, see section 4).

\

{\it Proof.} \hs Suppose $L$ has a CP-ideal $P$.  So, $P$ is a
5-dimensional abelian ideal of $L$.  Now take the linear
functional \hs $f = e_7^\ast \in L^\ast$, \hs which is regular.
Put \hs $A = ke_8 \subset Z(L)$. \hs This is a 1-dimensional ideal
of $L$ contained in $P$ and \hs $f(A) = 0$.  \hs By Lemma 1.13 \hs
$Q = P/A$ \hs is a CP-ideal of $L/A$.  Clearly, $L/A$ is a 7-
dimensional nilpotent Lie algebra of index 1, with basis \hs $x_1
= e_1 + A, \ldots, x_7 = e_7 + A$. \hs So, $Q$ is a 4-dimensional
abelian ideal of $L/A$. One verifies that there are \hs $\lambda,
\mu \in k$, \hs not both zero such that $Q$ is generated by \hs
$\lambda x_1 + \mu x_4, x_5, x_6, x_7$. \hs Then $P$ is generated
by \hs $\lambda e_1 + \mu e_4, e_5, e_6, e_7, e_8$.  \hs But this
contradicts the fact that $P$ is commutative, since
\begin{eqnarray*}
[\lambda e_1 + \mu e_4, e_5] = -\lambda e_8 \hs \mbox{and}\hs [\lambda e_1 + \mu
e_4, e_6] = -\mu e_8
\end{eqnarray*}

\item[(ii)] Let $V$ be a vector space over $k$ with basis \hs $e_1,\ldots, e_n$;
\hs $n \geq 2$. \hs Take the vector space $\bigwedge\limits^2 V$ with basis \hs
$e_{ij} = e_i \wedge e_j$, \hs $i < j$. \hs Next, consider the Lie algebra
\begin{eqnarray*}
L = V \oplus \ts\bigwedge\limits^2 V
\end{eqnarray*}
with nonvanishing brackets \hs $[e_i, e_j] = e_{ij}$, \hs $i < j$.  \hs Clearly,
\hs $[L,L] = \bigwedge\limits^2 V = Z(L)$. \hs So, $L$ is 2-step nilpotent of
dimension \hs $n + \frac{1}{2}n (n-1) = \frac{1}{2} n(n+1)$. \hs Let \hs $x,y
\in V$. \hs Then it is easy to see that
\begin{eqnarray*}
[x,y] = 0 \hs \Leftrightarrow \hs x,y\hs \mbox{are linearly dependent over
$k$}\insp (\ast)
\end{eqnarray*}
Next, we take $n$ to be even.  Then, $\rank_{R(L)} ([e_i, e_j]) = n$. \hs This
implies that
\begin{eqnarray*}
i(L) = \dim L - n = \frac{1}{2} n(n-1) = \dim Z(L)
\end{eqnarray*}
i.e. $L$ is square integrable.
\end{itemize}
Finally, take \hs $n = 4$. \hs Then \hs $\dim L = 10$, \hs $\dim
Z(L) = i(L) = 6$\hs and \\ $\frac{1}{2} (\dim L + i(L)) = 8.$\hs
But, because of $(\ast)$, $L$ has no 8-dimensional abelian Lie
subalgebra containing $Z(L)$, \hs i.e. $L$ has no CP's.  The same
holds for all even $n \geq 4$, using a similar argument.

\

{\bf THEOREM 3.2} \hs Let $L$ be a Lie algebra having a
nondegenerate, invariant bilinear form b.  Then \hs $F(L) = L$.
\hs In particular, $L$ does not admit a CP unless $L$ is abelian.

\

{\it Proof.} \hs Take \hs $y \in L$ \hs and consider the map
$\varphi_y$ sending each \hs $x \in L$ \hs into \hs $b(x,y)$.
Clearly, \hs $\varphi_y \in L^{\ast}$ \hs and the map \hs $\varphi
: L \rightarrow L^\ast$ \hs sending $y$ into $\varphi_y$ is an
isomorphism of $L$-modules. Consequently, $y$ and $\varphi_y$
have the same stabilizer in $L$, i.e. \hs $C(y) = L(\varphi_y)$.\\
Next, put \hs $\Omega = \varphi^{-1} (L^\ast_{\reg})$.\hs Then,
\begin{eqnarray*}
F(L) = \ts\sum\limits_{f\in L_{\reg}^\ast} L(f) = \ts\sum\limits_{y \in
\Omega}L(\varphi_y) = \ts\sum\limits_{y \in \Omega} C(y)
\end{eqnarray*}
Clearly, $F(L)$ contains $\Omega$, which is an open dense subset
of $L$ for the Zariski topology since $\varphi$ is a linear
isomorphism.  Consequently, \hs $F(L) = L$.

\

{\bf EXAMPLES 3.3}
\begin{itemize}
\item[(1)] $L$ semi-simple (take $b$ to be the Killing form of $L$).
\item[(2)] The diamond Lie algebra with basis \hs $t, x, y, z$ \hs and
nonvanishing brackets $[t,x] = -x$, \hs $[t,y] = y$ \hs and \hs $[x,y] = z$. \hs
Let $b$ be the symmetric bili\-near form with nonzero entries \hs $b(t,z) = 1$
\hs and \hs $b(x,y) = -1$.
\item[(3)] Let {\myfont g}$_5$ be the 5-dimensional nilpotent Lie algebra over
$k$ with basis \hs $x_1,\ldots, x_5$\hs and nonvanishing brackets \hs $[x_1,x_2]
= x_3, \hs [x_1, x_3] = x_4, \hs [x_2, x_3] = x_5$.\\
Let $b$ be the symmetric bilinear form with nonzero entries:
\begin{eqnarray*}
b(x_1, x_5) = b(x_3, x_3) = 1\hs \mbox{and}\hs b(x_2, x_4) = -1
\end{eqnarray*}
\item[(4)] Let {\myfont g}$_6$ be the 6-dimensional 2-step nilpotent Lie algebra
with basis \hs $x_1,\dots, x_6$ \hs and nonvanishing brackets \hs $[x_1,x_2] =
x_6, \hs [x_1,x_3] = x_4, \hs [x_2, x_3] = x_5$.\\
Let $b$ be the symmetric bilinear form with nonzero entries:
\begin{eqnarray*}
b(x_1,x_5) = b(x_3, x_6) = 1 \hs \mbox{and}\hs b(x_2, x_4) = -1
\end{eqnarray*}
(see [B1, p.133])
\item[(5)] Consider the semi-direct product \hs $L = sl(2,k) \oplus W_2$, \hs
where $W_2$ is the 3-dimensional irreducible $sl(2,k)$-module. \hs $L$ also
admits a nondegenerate, invariant, symmetric bilinear form.
\end{itemize}
\ \\
{\bf PROPOSITION 3.4} \hs Among all the different types of indecomposable
nilpotent Lie algebras over $\IC$ of dimension $n \leq 7$, only the following 30
Lie algebras do not have a CP:
\begin{itemize}
\item[1)] $n=5$: \hs {\myfont g}$_5$ (see (3) of Examples 3.3)
\item[2)] $n=6$: \hs From Morozov's classification [Mo, p.168] the Lie algebras
\hs $3(\cong${\myfont g}$_6)$, 21 and 22.
\item[3)] $n=7$: \hs From Seeley's classification [See]: 2, 5, $7_K$; \hs 2, 5,
$7_L$; \hs 2, 4, $7_D$; \\
2, 4, $7_E$;\hs 2, 4, $7_G$; \hs 2, 4, $7_H$; \hs 2, 4, $7_J$; \hs 2, 4, $7_K$;
\hs 2, 4, $7_Q$; \hs 2, 4, $7_R$;\\
2, 3, 5, $7_C$; \hs 2, 3, 5, $7_D$; \hs 2, 3, 4, 5, $7_B$; \hs 2, 3, 4, 5,
$7_C$; \hs 2, 3, 4, 5, $7_D$;\\
2, 3, 4, 5, $7_F$; \hs 2, 3, 4, 5, $7_G$; \hs 1, 3, 5, $7_S$, \hs $\xi = 1$; \hs
1, 3, 4, 5, $7_H$; \\
1, 2, 4, 5, $7_C$; \hs 1, 2, 4, 5, $7_F$; \hs 1, 2, 4, 5, $7_H$; \hs
1, 2, 4, 5, $7_K$; \hs 1, 2, 4, 5, $7_L$; \\
1, 2, 4, 5, $7_N$, \hs $\xi = 1$; \hs 1, 2, 3, 4, 5, $7_I$, \hs $\xi = 0$\\
Note that the infinite families fail to have a CP only for exceptional values of
the parameter $\xi$.
\end{itemize}

{\it Proof:} \hs This is done case by case, considering only the
ones that are not square integrable (Proposition 2.2).  Usually,
CP's are easy to spot by loo\-king at the multiplication table. To
prove that a Lie algebra $L$ has no CP's is more difficult
however.  This can be achieved by using Proposition 1.9 or by
showing that $F(L)$ is not commutative.  For instance, take $L =
\hs 1, 2, 4, 5, 7_N$, $\xi = 1$. See Remark 2.3 for its Lie
brackets.  One verifies that $F(L) = \langle a-b, c, d, e, f,
g\rangle$, which is not commutative.

\

{\bf REMARK 3.5} \hs Having a CP is not preserved under
degeneration (for a definition we refer to [GO1] or [GO2]).
Indeed, $g_5$, which has no CP's (see 3 of Examples 3.3), is a
degeneration of the Lie algebra $h_5$ with basis \hs $x_1,\ldots,
x_5$\hs over $\IC$ and nonzero brackets \hs $[x_1,x_2] = x_3$, \hs
$[x_1,x_3] = x_4$, \hs $[x_1, x_4] = x_5$ \hs and \hs $[x_2, x_3]
= x_5$ \hs for which \hs $\langle x_3, x_4, x_5\rangle$ \hs is a
CP.  On the other hand, the Lie algebra $j_5$ with the same basis
and nonzero brackets \hs $[x_1,x_2] = x_3$ \hs and \hs $[x_1, x_3]
= x_4$ \hs admits a CP (namely \hs $\langle x_2, x_3, x_4,
x_5\rangle$) and is a degeneration of {\myfont g}$_5$ \hs [GO1,
p.323].

\

\begin{center}
{\bf 4. CP-IDEALS}
\end{center}
These are by far the most interesting CP's.  The following shows
that they occur as often as ordinary CP's, at least in the
solvable case.

\

{\bf THEOREM 4.1} \hs Let $L$ be solvable and $k$ algebraically
closed.  Let $m$ be the maximum dimension of all abelian ideals of
$L$.  Clearly, \hs $m \leq \frac{1}{2} (\dim L +i(L))$ \hs [O2, p.
706].  Then the following are equivalent:
\begin{itemize}
\item[(1)] $L$ admits a CP
\item[(2)] $L$ admits a CP-ideal
\item[(3)] $m = \frac{1}{2} (\dim L + i(L))$
\end{itemize}

{\it Proof.} \hs It suffices to show that \hs (1) $\Rightarrow$
(2), \hs since \hs (2) $\Rightarrow$ (1) \hs and \hs (2)
$\Leftrightarrow$ (3) \hs are clear. Let $G$ be the adjoint
algebraic group of $L$, i.e. the smal\-lest algebraic subgroup of
\hs $\Aut L$ \hs such that $L(G)$ contains $\ad L$ [D, 1.1.14].
Clearly, $\ad L$ and hence its algebraic hull $L(G)$ are solvable
(since they have the same derived algebra [Ch, p.173], which is
nilpotent).  Therefore $G$ is a solvable connected group.  Next
put $p = \frac{1}{2}(\dim L + i(L))$. \hs Then the set $C$ of all
CP's is a nonempty (by assumption) closed subset of the
Grassmannian \hs $\Gr(L,p)$, \hs which is an irreducible and
complete algebraic variety [D, 1.11.8-9].  Hence $C$ is also
complete.  Now $G$ acts morphically on $C$, mapping each CP $H$ on
$g(H)$, \hs $g \in G$. \hs By Borel's theorem, $G$ has a fixed
point $P$ in $C$ [Bo, p.242].  So, \hs $g(P) = P$\hs for all $g
\in G$. In particular, $\ad x (P) \subset P$ \hs for all $x \in
L$.  Consequently, $P$ is a CP-ideal of $L$.

\

{\bf REMARK 4.2} \hs (a) The number $m$ is an important
characteristic of a Lie
algebra, often used in classifications.\\
(b) It is now easy to see that the 8-dimensional Lie algebra (i)
of 3.1 has no CP's (go over to the algebraic closure of $k$ and
use Theorem 4.1).

\

{\bf THEOREM 4.3} \hs Let $P$ be an ideal of a Lie algebra $L$ and
let $P$ be a pola\-rization of $L$ with respect to some \hs $f \in
L^\ast$.  \hs Then we have
\begin{itemize}
\item[(1)] If \hs $f \in L_{\reg}^\ast$ \hs then $P$ is solvable (in fact \hs
$P'' = 0$).  If in addition $L$ is Frobenius or nilpotent of index one, then $P$
is a CP-ideal of $L$.
\item[(2)] If $k$ is algebraically closed and \hs $f \in L_{\reg}^\ast$, \hs
then the induced representation \hs $\ind(f|_P,L)$ \hs is simple.
\item[(3)] If $L$ is completely solvable then $P$ is a Vergne polarization.  In
particular, \hs $\ind(f|_P,L)$ \hs is absolutely simple.
\end{itemize}

{\it Proof.}
\begin{itemize}
\item[(1)] Take \hs $x \in L$ \hs and \hs $y, y' \in P$, \hs then
\begin{eqnarray*}
f([x,[y,y']]) &=& f([[x,y],y']) + f([y,[x,y']])\\
&=& 0 \hs \mbox{since $P$ is an ideal}
\end{eqnarray*}
and \hs $f([P,P]) = 0$.

\

Hence, \hs $[y,y'] \in L(f)$. \hs Therefore, \hs $P' =
[P,P]\subset L(f)$. \hs This implies that \hs $P'' = 0$ \hs since
$L(f)$ is abelian by [D,1.11.7].  Now, suppose $L$ is Frobenius,
i.e. \hs $i(L) = 0$. \hs Then \hs $L(f) = 0$ \hs which forces \hs
$[P,P] = 0$. \hs On the other hand, if $L$ is nilpotent of index
1, then \hs $\dim L(f) = 1$.  \hs We may assume that $f \neq 0$.
Clearly, \hs $[P,P] \neq L(f)$ \hs since \hs $f([P,P]) = 0$ \hs
and \hs $f(L(f)) \neq 0$ \hs [BC, p.89].  So, we conclude that \hs
$[P,P] = 0$. \item[(2)] By [RV, p.395] or [D, 10.5.7] there exists
a solvable polarization $H$ of $L$ w.r.t. $f$ such that \hs $H
\cap P$ \hs is a solvable polarization of $P$ w.r.t. $f|_P$ and
such that the twisted induced representation \hs $\ind^{\sim}
(f|_H,L)$ \hs is simple.  First we observe that
\begin{eqnarray*}
\dim H = \ts\frac{1}{2} (\dim L + \dim L(f)) = \dim P \insp \insp (\bullet)
\end{eqnarray*}
Similarly,
\begin{eqnarray*}
\dim (H \cap P) = \ts\frac{1}{2} (\dim P + \dim P(f|_P)) = \dim P
\end{eqnarray*}
since \hs $P(f|_P) = \{x \in P \mid f([x,P]) = 0\} = P$. \hs It follows that \hs
$H \cap P = P$, \hs i.e. \hs $P \subset H$. \hs Hence, by $(\bullet)$, we see
that \hs $P = H$.\\
Consequently, \hs $\ind^\sim(f|_P, L)$\hs is simple.  Finally, \hs
$\ind^\sim(f|_P, L) = \ind(f|_P,L)$ \hs because $P$ is an ideal of $L$ [D,
5.2.1].
\item[(3)] $L$ being completely solvable, we can find a flag of ideals of $L$:
\begin{eqnarray*}
L = L_n \supset \ldots \supset L_p \supset \ldots \supset L_1 \supset L_0 = (0)
\end{eqnarray*}
such that \hs $L_p = P$ where $p = \dim P$. Put $f_i = f|_{L_i}$ and $P_j =
\sum\limits_{i \leq j} L_i(f_i)$.\\
Then $P_n$ is the so called Vergne polarization w.r.t. this flag and $f \in
L^\ast$ [BGR, 9.4].  We claim that \hs $P = P_n$. \hs Clearly,
\begin{eqnarray*}
L_i(f_i) = \{x \in L_i \mid f([x,L_i]) = 0\} = L_i \cap L_i^\perp
\end{eqnarray*}
In particular, \hs $L_p (f_p) = L_p \cap L_p^\perp = P \cap P^\perp = P$ since
\hs $P = P^\perp$ \hs w.r.t. \hs $f \in L^\ast$. \hs This implies that \hs $P
\subset P_n$. \hs On the other hand consider \hs $L_j (f_j)$.\\
If $j \leq p$, then \hs $L_j(f_j) \subset L_j \subset L_p = P$.\\
If $j > p$, then \hs $P = L_p \subset L_j$ \hs implies that \hs $L_j(f_j) = L_j
\cap L_j^\perp \subset L_j^\perp \subset P^\perp = P$. \hs Consequently, \hs
$P_n = \sum\limits_{j=1}^n L_j(f_j) \subset P$.
\end{itemize}
\ \\
{\bf COROLLARY 4.4} \hs Let $P$ be a CP-ideal of a Lie algebra $L$ and take
any\\
$f \in L_{\reg}^\ast$.  Then,
\begin{itemize}
\item[1.] If $k$ is algebraically closed, then \hs $\ind(f|_{P}, L)$ \hs is
simple.
\item[2.] If $L$ is completely solvable, then $P$ is a Vergne polarization
w.r.t. $f$ and any flag of ideals containing $P$.  In particular, $\ind(f|_P,
L)$ is absolutely simple.
\item[3.]
\begin{itemize}
\item[(a)] $Sz(U(L)) \subset U(P)$ \hs and \hs $Sz(D(L)) \subset D(P)$ \hs where
\hs $Sz(U(L)) = \bigoplus\limits_\lambda U(L)_\lambda$ \hs is the semi-center of
$U(L)$.  Similarly for \hs $Sz(D(L))$. \hs This generalizes [D, 6.1.6].
\item[(b)] Put \hs $\wedge (L) = \{\lambda \in L^\ast\hs \mid \hs U(L)_\lambda
\neq 0\}$ \hs and \hs $L_\wedge = \bigcap\limits_{\lambda \in \wedge(L)} \ker
\lambda$. \hs Then, \hs $P \subset L_\wedge$.
\end{itemize}
\end{itemize}

\

{\it Proof.} \hs (1) and (2) follow directly from Theorem 4.3.\\
(3) Let \hs $u \in U(L)_\lambda$ \hs be any semi-invariant with weight \hs
$\lambda \in \wedge (L)$, \hs i.e. \hs $[x,u] = \lambda(x)u$ \hs for all \hs $x
\in L$.\\
Now, take \hs $x \in P$. \hs Then \hs $adx(L) \subset P$ \hs and
\hs $(adx)^2 = 0$ \hs since $P$ is a commutative ideal of $L$. So,
\hs $adx$ \hs is nilpotent. This implies that  \hs $\lambda (x) =
0$ \hs and \hs $[x,u] = 0$.\hs Consequently, \hs $x \in L_\wedge$
\hs which shows (b) and also \hs $u \in C(U(P)) = U(P)$. \hs
Therefore, \hs $Sz(U(L)) \subset U(P)$. \hs Similarly for \hs
$Sz(D(L)) \subset D(P)$ \hs (since \hs $C(D(P)) = D(P))$.

\

 {\bf REMARK 4.5} \hs  The previous corollary does not hold for
arbitrary CP's of $L$.  For example, let $L$ be the 2-dimensional
Lie algebra over an algebraically closed field $k$ with basis
$x,y$ and nonzero bracket \hs $[x,y] = y$. \hs $L$ is Frobenius
and \hs $f \in L^\ast$ \hs with \hs $f(x) = 0$ \hs and $\hs f(y) =
1$ \hs is regular.  Clearly, \hs $P = kx$ \hs is a CP of $L$
w.r.t. \hs $f \in L^\ast$. \hs But \hs $\ind^\sim(f|_P, L)$\hs is
not simple [BGR, p.95].  Also, $y$ is a semi-invariant of $L$ but
\hs $y \notin U(P)$.

\

 The following, which we recall from [O2, p.708], describes how
CP-ideals naturally arise in certain semi-direct products.

\

{\bf PROPOSITION 4.6.} \hs Let {\myfont g} be a Lie algebra with basis \hs
$\{x_1,\ldots, x_m\}$ \hs and let $V$ be a {\myfont g}-module with basis \hs
$\{v_1,\dots, v_n\}$ \hs with \hs $\dim${\myfont g}\ $\leq \dim V$. \hs For each
\hs $f \in V^\ast$ \hs we put
\begin{eqnarray*}
\mbox{{\myfont g}}(f) = \{x \in \mbox{{\myfont g}}\hs \mid \hs f(xv) = 0 \hs
\mbox{for all} \hs v \in V\}
\end{eqnarray*}
the stabilizer of $f$.  Consider the semi-direct product \hs $L =$\ {\myfont
g}$\ \oplus\ V$ \hs in which \hs $[x,v] = xv$, \hs $x \in$\ {\myfont g}, \hs $v
\in V$ \hs and in which $V$ is an abelian ideal.  Then the following are
equivalent:
\begin{itemize}
\item[(1)] $D(V) \hs (=R(V))$ is a maximal subfield of $D(L)$
\item[(2)] $V$ is a CP-ideal of $L$
\item[(3)] $i(L) = \dim V - \dim$\ {\myfont g}
\item[(4)] $\rank_{R(V)} (e_iv_j) = \dim$\ {\myfont g}
\item[(5)] {\myfont g}$(f) = 0$ \hs for some \hs $f \in V^\ast$
\end{itemize}

\

{\bf REMARK 4.7} \hs If $k$ is algebraically closed, {\myfont g} a
simple Lie algebra, acting irreducibly on $V$, then the conditions
of the proposition are satisfied if and only if \\ $\dim$\
{\myfont g} \ $< \dim V$. \hs [AVE, p.196].

\

 The following shows that if a Lie algebra $L$ admits a CP-ideal
then its structure comes close to that of the semi-direct product
considered in Proposition 4.6.

\

 {\bf COROLLARY 4.8} \hs Let $V$ be a commutative ideal of $L$.
Clearly, the Lie algebra \hs {\myfont g} $= L/V$ \hs acts on $V$.
Consider the semi-direct product \hs $L_1 =$\ {\myfont g}$\ \oplus
\ V$. \hs Then,
\begin{center}
$V$ is a CP of L \hs $\Leftrightarrow \hs V$ is a CP of $L_1$
\end{center}
In that case, \hs $i(L_1) = i(L)$.

\

 {\it Proof.} \hs Let \hs $g \in L^\ast$ \hs and put \hs $f =
g|_V \in V^\ast$. \hs Then, we claim that \hs {\myfont g}$(f) =
V^g/V$. \hs Indeed,
\begin{eqnarray*}
\overline{x} = x + V \in \ \mbox{{\myfont g}}(f) \hs &\Leftrightarrow& \hs
f([\overline{x}, V]) = 0\\
&\Leftrightarrow& \hs f([x,V]) = 0\\
&\Leftrightarrow& \hs g([x,V]) = 0\\
&\Leftrightarrow& \hs x \in V^g \hs \Leftrightarrow \hs \overline{x} \in V^g/V
\end{eqnarray*}
We now proceed with the proof\\
$\Rightarrow :$ \hs $\dim V = \frac{1}{2} (\dim L + i(L))$. \hs Also, \hs $V^g =
V$ \hs for some \hs $g \in L^\ast$ \hs by Lemma 1.1.  Hence, \hs {\myfont g}$(f)
= 0$. \hs By Proposition 4.6 \hs $V$ is a CP of $L_1$ and
\begin{eqnarray*}
i(L_1) = \dim V - \dim \ \mbox{{\myfont g}} &=& \dim V - (\dim L - \dim V)\\
&=& 2\dim V - \dim L = i(L).
\end{eqnarray*}
$\Leftarrow :$ By Proposition 4.6, \hs {\myfont g}$(f) = 0$ \hs
for some \hs $f \in V^\ast$. \hs Next, choose \hs $g \in L^\ast$
\hs such that \hs $f = g|_V$. \hs Then, \hs $V^g/V =$ {\myfont
g}$(f) = 0$.  \hs So, \hs $V^g = V$ \hs which by Lemma 1.2 implies
that $V$ is a CP of $L$.

\

\begin{center}
{\bf 5. CP-IDEALS IN CERTAIN FROBENIUS LIE ALGEBRAS}
\end{center}
Let $L$ be a Frobenius Lie algebra with a CP-ideal $P$.  Take any
\hs $f \in L^\ast_{reg}$ \hs and assume that $k$ is algebraically
closed.  Then \hs $I(f) = 0$ \hs by [O1, p.42].  So, by Corollary
4.4 \hs $\ind(f|_P,L)$ \hs is a faithful irreducible
representation of $U(L)$.  Next, let \hs $x_1, \ldots, x_m$,
$y_1,\ldots, y_m$ \hs be a basis of $L$ such that \hs $y_1,\ldots,
y_m$ is a basis of $P$.  Then \hs $\det ([x_i, y_j]) \in S(P)$ \hs
is a nonzero semi- invariant under the action of \hs $\Aut L$ \hs
[O1,p.28].  It is also known that Frobenius Lie algebras give rise
to constant solutions for the classical Yang- Baxter equation
[BD].

\

The following is a special case of Proposition 4.6.

 \

 {\bf
COROLLARY 5.1} \hs Let {\myfont g} be a Lie algebra and $V$ a
{\myfont g}-module such that \hs $\dim$\ {\myfont g} $= \dim V$.
\hs Consider the semi- direct product \hs $L =$ {\myfont g}$\
\oplus \ V$.  Then the following are equivalent:
\begin{itemize}
\item[(1)] $R(V)$ is a maximal subfield of $D(L)$
\item[(2)] $V$ is a CP-ideal of $L$
\item[(3)] $L$ is Frobenius
\item[(4)] {\myfont g}$(f) = 0$ \hs for some \hs $f \in V^\ast$
\end{itemize}

{\bf EXAMPLE 5.2} \hs Let {\myfont g} be Frobenius and let \hs $V
=$ {\myfont g} be the adjoint representation.

\

 {\bf EXAMPLE 5.3} \hs The above condition is satisfied if
{\myfont g} is reductive over an algebraically closed field $k$
and $V^\ast$ is a prehomogeneous {\myfont g}-module (i.e. $V^\ast$
has an open {\myfont g}-orbit) with $\dim$ {\myfont g} $= \dim V$.
These modules have been studied extensively by the Japanese school
since 1977 [SK], [KKTI].

\

 {\bf EXAMPLE 5.4} \hs Let $A$ be a left-symmetric algebra
(LSA), i.e. a finite dimensional vector space provided with a
bilinear product \hs $A \times A \rightarrow A$, $(a,b)
\rightarrow ab$ \hs which satisfies
\begin{eqnarray*}
a(bc) - (ab)c = b(ac) - (ba) c \insp \insp \insp (\ast)
\end{eqnarray*}
for all \hs $a,b,c \in A$. \hs  There is an extensive literature on LSA's, see
for example [H], [Seg]. Vinberg used LSA's to classify convex homogeneous cones
[V].  A left-symmetric algebra is Lie-admissable.  This means that $A$ becomes a
Lie algebra, which we denote by {\myfont g}, for the Lie bracket \hs $[a,b] = ab
- ba$, \hs $a,b \in A$. \hs Using $(\ast)$ we observe that
\begin{eqnarray*}
[a,b]c = (ab)c - (ba)c = a(bc) - b(ac).
\end{eqnarray*}
Therefore, $A$ becomes a {\myfont g}-module, which we denote by $V$, for the
bilinear map
\begin{eqnarray*}
\mbox{{\myfont g}} \times V \rightarrow V, \hs (x,v) \rightarrow xv
\end{eqnarray*}
Now, suppose $A$ contains a nonzero element \hs $f \in A$ \hs which is not a
right zero divisor of $A$.  Let $V^\ast$ be the dual module of $V$.  Identifying
the module $V^{\ast\ast}$ with $V$, we may consider $f$ to be an element of
$(V^{\ast})^\ast$. Clearly, the stabilizer \hs {\myfont g}$(f) = \{x \in
\mbox{{\myfont g}}\ \mid xf = 0\} = 0$ \hs by assumption.\\
Finally, using Corollary 5.1 we may conclude that the semi-direct product \\
$L =$ {\myfont g} \ $\oplus \ V^\ast$ \hs is a Frobenius Lie
algebra in which $V^\ast$ is a CP-ideal.

\

 {\bf REMARK 5.5} \hs In characteristic $p > 2$ a similar result
can be obtained if $A$ is a finite dimensional simple Novikov
algebra and where $V$ is a certain irreducible $A$-module. We
recall that a nonassociative $k$-algebra is said to be a left
Novikov algebra if $A$ is left symmetric, satisfying the identity
\hs $(ab)c = (ac)b$ \hs for all \hs $a,b,c \in A$. \hs In
characteristic zero E. Zelmanov showed that finite dimensional
simple Novikov algebras are all one- dimensional [Z].  Recently
simple Novikov algebras and their irreducible modules have been
determined by M. Osborn and X. Xu [Os], [X].

\

We now focus on a special case, which provides an interesting link
between Frobenius algebras and Frobenius Lie algebras.

\

 {\bf PROPOSITION 5.6} \hs Let $A$ be a finite dimensional
associative algebra over $k$ with a unit element.  A becomes a Lie
algebra {\myfont g} for the Lie bracket \hs $[a,b] = ab - ba$, \hs
$a,b \in A$, \hs and \hs $V = A$ \hs becomes a {\myfont g}-module
by left multiplication.  Consider the semi-direct product \hs $L
=$ {\myfont g}\ $\oplus \ V$. \hs Then the following conditions
are equivalent:
\begin{itemize}
\item[(1)] $A$ is a Frobenius algebra
\item[(2)] $L$ is a Frobenius Lie algebra
\item[(3)] $V$ is a CP-ideal of $L$
\item[(4)] $R(V)$ is a maximal subfield of $D(L)$
\end{itemize}
{\it Proof.} \hs In view of Corollary 5.1 it suffices to show that (1) is
equivalent with \hs {\myfont g}$(f) = 0$ \hs for some $f \in V^\ast$. So, take
\hs $f \in V^\ast$. \hs Then
\begin{eqnarray*}
\mbox{{\myfont g}}(f) = \{a \in A \mid f(ab) = 0 \ \mbox{for all} \ b \in A\}
\end{eqnarray*}
Clearly, \hs {\myfont g}$(f) = 0$ \hs if and only if the bilinear
map \hs $A \times A \rightarrow k$, \hs $(a,b) \rightarrow f(ab)$
\hs is nondegenerate, i.e. $A$ is a Frobenius algebra [CR, Theorem
61.3].

\

 Finally, we devote our attention to certain Frobenius Lie
subalgebras of a semi- simple Lie algebra.

\

{\bf THEOREM 5.7} \hs Let $L$ be a semi-simple Lie algebra of rank
$r$ over $k$, $k$ algebraically closed, and let $x$ be a principal
nilpotent element of $L$ (i.e. the centralizer $C(x)$ of $x$ in
$L$ has dimension $r$).  Then the normalizer $F$ of $C(x)$ in $L$
is a solvable Frobenius Lie subalgebra of $L$ in which $C(x)$ is a
CP-ideal.

\

 {\it Proof.} \hs It is well known that $C(x)$ is abelian [K].
Clearly, $C(x)$ is an ideal of $F$.  In 1991 R. Brylinski and B.
Kostant showed that \hs $\dim F = 2r$ \hs and that $F/C(x)$, \hs
and hence also $F$, is solvable $[BK]$.
Recently, D. Panyushev proved that $F$ is Frobenius [P2, Theorem 5.5].

\

\begin{center}
{\bf 6. CP-IDEALS IN THE NILRADICAL OF PARABOLIC LIE SUBALGEBRAS OF A SIMPLE LIE
ALGEBRA}
\end{center}
{\bf THEOREM 6.1} \hs Let $B$ be a Borel subalgebra of a simple Lie algebra $L$
over $k$, $k$ algebraically closed, of rank $r$ and let $N$ be the nilradical of
$B$.  Then,
\begin{itemize}
\item[(1)] $N$ admits a CP \hs $\Leftrightarrow$ \hs $L$ is of type $A_r$ or
$C_r$.\\
In these 2 cases $N$ has a CP-ideal $P$, which is an ideal of $B$.
\item[(2)] $P$ is also a CP-ideal of $B$ in case $L$ is of type $C_r$, \hs $r
\geq 1$.
\end{itemize}

{\it Proof.} \hs The information on \hs $i(N)$, \hs $i(B)$ \hs in table 1 is
obtained from [E1], [E2].  Also, we know that \hs $i(N) + i(B) = r$ \hs [P2,
1.5].

\begin{center}
\begin{tabular}{|cc|c|c|c|c|c|}
\hline
 & &$\dim N$ &$i(N)$ &$i(B)$ &$\frac{1}{2}(\dim N + i(N))$ &$m$\\
\hline
$A_{2t}$   &$t \geq 1$ &$t(2t+1)$ &$t$ &$t$ &$t(t+1)$ &$t(t+1)$\\
$A_{2t+1}$ &$t \geq 0$ &$(t+1)(2t+1)$ &$t+1$ &$t$ &$(t+1)^2$ &$(t+1)^2$\\
\hline
$B_3$ & &9 &3 &0 &6 &5\\
$B_r$ &$r \geq 4$ &$r^2$ &$r$ &0 &$\frac{1}{2}r(r+1)$ &$\frac{1}{2}r(r-1) + 1$\\
\hline
$C_r$ &$r \geq 2$ &$r^2$ &$r$ &0 &$\frac{1}{2} r(r+1)$ &$\frac{1}{2}r(r+1)$\\
\hline
$D_{2t}$   &$t \geq 2$ &$2t(2t-1)$ &$2t$ &0 &$2t^2$ &$t(2t-1)$\\
$D_{2t+1}$ &$t \geq 2$ &$2t(2t + 1)$ &$2t$ &1 &$2t(t+1)$ &$t(2t+1)$\\
\hline
\multicolumn{2}{|c|}{$E_6$} &36 &4 &2 &20 &16\\
\hline
\multicolumn{2}{|c|}{$E_7$} &63 &7 &0 &35 &27\\
\hline
\multicolumn{2}{|c|} {$E_8$} &120 &8 &0 &64 &36\\
\hline
\multicolumn{2}{|c|} {$F_4$} &24 &4 &0 &14 &9\\
\hline
\multicolumn{2}{|c|} {$G_2$} &6 &2 &0 &4 &3\\
\hline
\end{tabular}\\
\ \\
\ \\
{\bf Table 1}
\end{center}
The idea is to compare the maximum dimension $m$ of abelian Lie subalgebras of
$N$, computed by Malcev [Ma, p.216] with the number \hs $\frac{1}{2}(\dim N +
i(N))$. \hs Then $N$ contains a CP if and only if these numbers coincide.
According to the table this occurs precisely if $L$ is of type $A_r$ or $C_r$.\\
Furthermore, we know from [PR, Table 1] that in both types ($A_r$ or $C_r$) $B$
has a maximal abelian ideal $P$ of dimension \hs $\frac{1}{2} (\dim N + i(N))$.
\hs Clearly $P \subset N$.  Therefore $P$ is a CP-ideal of $N$.  This can also
be deduced from Theorem 4.1.\\
(2) Using Lemma 1.3 we see that $P$ is also a CP-ideal of $B$ if and only if
\begin{eqnarray*}
i(N) &=& i(B) +  \dim B - \dim N\\
&\Leftrightarrow& \hs i(N) - i(B) = r\\
&\Leftrightarrow& \hs i(B) = 0 \hs \mbox{(since \hs $i(N) + i(B) = r$).}
\end{eqnarray*}
and this happens when $L$ is of type \hs $A_1 (=C_1)$ or $C_r$, \hs $r \geq
2$.

\

{\bf THEOREM 6.2} \hs Let $L$ be a simple Lie algebra over $k$, $k$
algebraically closed, of type $A_r$ or $C_r$, $\pi$ a parabolic Lie subalgebra
of $L$.  Then the nilradical $N$ of $\pi$ admits a CP-ideal $P$.  Furthermore,
\begin{itemize}
\item[(1)] suppose $L$ is of type $A_r$ and $\pi$ of type \hs
$(p_1,\ldots,p_m)$. \hs Put \hs $n=r+1$ \hs and \hs $p = p_1 + \ldots + p_\ell$,
\hs $1 \leq \ell \leq m$, \hs such that \hs $\left|\sum\limits_{i=1}^\ell p_i -
\ds\frac{n}{2}\right|$ \hs is as small as possible.  Then,
\begin{eqnarray*}
i(N) = 2p (n-p) - \frac{1}{2} \left(n^2 - \sum\limits_{i=1}^m p_i^2\right)
\end{eqnarray*}
\item[(2)] suppose $L$ is of type \hs $C_r$, $r \geq 2$, \hs and $\pi$ of type
\hs $(p_1,\ldots, p_m)$. \hs Put \hs $\ell = \left[\frac{m}{2}\right]$, \hs then
\begin{eqnarray*}
i(N) = \frac{1}{2} \sum\limits_{i=1}^\ell p_i (p_i + 1)
\end{eqnarray*}
\end{itemize}

{\bf REMARK 6.3}
\begin{itemize}
\item[a)] The first formula is new.  A recursive formula for $i(N)$ was already
established in [E1].  A different proof for the second formula can also be found
in [E1].
\item[b)] (made by the referee) A. Joseph already gave a formula for $i(N)$ in
an arbitrary simple Lie algebra, using a maximal subset of strongly orthogonal
positive roots [J, (ii) of Proposition 2.6].  Being applied to $A_r$ or $C_r$,
Joseph's formula gives the above explicit expressions.
\end{itemize}
{\it Proof.} \hs (1) Let \hs $L = \mbox{sl}(V)$ \hs where $V$ is an $n$-
dimensional vector space over $k$.  By [B2, p.187] we can find a flag $F$ of
subspaces of $V$:
\begin{eqnarray*}
\{0\} = F_0 \subset F_1 \subset \ldots \subset F_m = V, \hs F_{i-1}
\subs\limits_{\neq} F_i
\end{eqnarray*}
such that $\pi$ (respectively its nilradical $N$) consists of all endomorphisms
\hs $x \in L$ \hs such that $xF_i \subset F_i$ \hs (resp. \hs $xF_i \subset
F_{i-1}$) \hs for \hs $1 \leq i \leq m$. \hs Put \hs $p_i = \dim (F_i/F_{i-1})$
\hs then \hs $\pi$ is said to be of type \hs $(p_1,\ldots, p_m)$.  \hs Next,
choose a basis \hs $e_1,\dots, e_n$ \hs of $V$ compatible with the flag $F$
(i.e. \hs $e_1,\ldots, e_{p_1} \in F_1\backslash F_0$, \hs etc.).  Then, $N$ can
be considered to be the Lie algebra of matrices of the form as shown in figure
1.
\begin{center}
\includegraphics{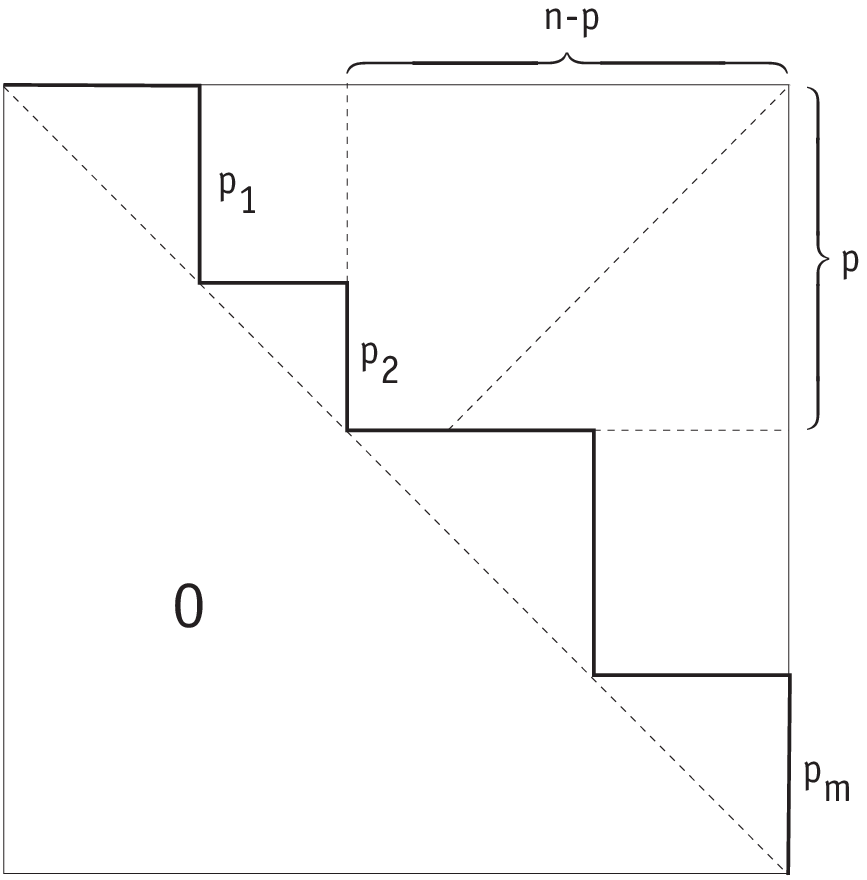}
\end{center}
\begin{center}
{\bf Figure 1}
\end{center}

We may assume, as is the case in figure 1, that \hs $p \leq \frac{n}{2}$ \hs
$(\ast)$.  In particular, \hs $p + p_{\ell +1} > \frac{n}{2}$. \hs As usual we
denote by $E_{ij}$ the $n \times n$ matrix whose $ij$-th entry is 1 and other
entries are zero.  Let $P$ be the subspace of $N$ generated by all $E_{ij}$ with
\hs $1 \leq i \leq p$; \hs $p+1 \leq j \leq n$. \hs So, $P$ consists of matrices
of the form \hs $\left(\begin{array}{cc}0 &M\\0 &0\end{array}\right)$ \hs where
$M$ is any \hs $p \times (n-p)$ \hs matrix.  It is easy to see that $P$ is an
abelian ideal of $N$.  We claim that $P$ is a CP of $N$.  Let \hs $f \in N^\ast$
\hs be defined by \hs $f(E_{p,n-p+1}) = \ldots = f(E_{1n}) = 1$ \hs and zero on
all other $E_{ij}$.  We want to show that \hs $P^f = P$.  Therefore we take \hs
$x \in P^f$. \hs We write
\begin{eqnarray*}
x = \ts\sum\limits_{i < j} \lambda_{ij} E_{ij} + y
\end{eqnarray*}
where \hs $E_{ij} \in N \backslash P$, \hs $\lambda_{ij} \in k$ \hs and \hs $y
\in P$. \hs We need to demonstrate that each \hs $\lambda_{i_0j_0} = 0$. \hs
There are two cases to distinguish:
\begin{itemize}
\item[(i)] $j_0 \leq p$.  \hs Then \hs $i_0 < j_0 \leq p$ \hs and \hs $s = (n+1)
- i_0 > (n+1) -p > p$. \hs Hence, \hs $E_{j_0s} \in P$ \hs and
\begin{eqnarray*}
0 &=& f([x, E_{j_0s}]) = \ts\sum\limits_{i < j} \lambda_{ij} f([E_{ij},
E_{j_0s}]) + f([y,E_{j_0s}])\\
&=& \ts\sum\limits_{i < j} \lambda_{ij} f(\delta_{jj_0} E_{is} - \delta_{si}
E_{j_0j})\\
&=& \ts\sum\limits_{i < j_0} \lambda_{ij_0} f(E_{is}) - \ts\sum\limits_{j > s}
\lambda_{sj} f(E_{j_0j})\\
&=& \lambda_{i_0j_0}
\end{eqnarray*}
($f(E_{j_0j}) = 0$ \hs since \hs $j_0 + j > i_0 + s = n + 1$)
\item[(ii)] $i_0 > p$ \hs and \hs $j_0 > p_1 + \ldots + p_\ell + p_{\ell+1} >
\ds\frac{n}{2}$.\\
By definition of $p$:
\begin{eqnarray*}
(p_1 + \ldots + p_\ell + p_{\ell +1}) - \frac{n}{2} \geq \ds\frac{n}{2} - p
\end{eqnarray*}
Hence
\begin{eqnarray*}
j_0 \geq (p_1 + \ldots + p_\ell + p_{\ell +1}) + 1 \geq n-p + 1
\end{eqnarray*}
So, \hs $t = (n+1) - j_0 \leq p < i_0$ \hs and \hs $E_{ti_0} \in P$.  Therefore
\begin{eqnarray*}
0 &=& f([E_{ti_0}, x]) = \ts\sum\limits_{i < j} \lambda_{ij} f([E_{ti_0},
E_{ij}]) + f([E_{ti_0}, y])\\
&=& \ts\sum\limits_{i < j} \lambda_{ij} f(\delta_{i_0i} E_{tj} - \delta_{jt}
E_{ii_0})\\
&=& \ts\sum\limits_{j > i_0} \lambda_{i_0j} f(E_{tj}) - \ts\sum\limits_{i < t}
\lambda_{it} f(E_{ii_0})\\
&=& \lambda_{i_0j_0}
\end{eqnarray*}
($f(E_{ii_0}) = 0$ \hs since \hs $i + i_0 < t + j_0 = n+1$).
\end{itemize}
In both cases: \hs $x = y \in P$.  \hs So, \hs $P^f \subset P$.
\hs Consequently, \hs $P^f = P$ \hs as the other inclusion is
obvious by the commutativity of $P$.

By Lemma 1.2 we may conclude that $P$ is a CP of $N$ and \hs $f \in
N_{\reg}^\ast$.\\
Finally, from \hs $\dim P = \frac{1}{2} (\dim N + i(N))$ \hs we obtain:
\begin{eqnarray*}
i(N) &=& 2 \dim P - \dim N\\
&=& 2p (n-p) - \frac{1}{2} (n^2 - \ts\sum\limits_{i=1}^m p_i^2)
\end{eqnarray*}

(2) Let \hs $L = \mbox{sp}(V)$ \hs where $V$ is a vector space over $k$ of
dimension $n = 2r$ provided with a nondegenerate alternating bilinear form
$\varphi : V \times V \rightarrow k$.  There exists an isotropic flag
\begin{eqnarray*}
\{0\} = F_0 \subset F_1 \subset \ldots \subset F_m = V
\end{eqnarray*}
i.e. \hs $F_i^\perp = F_{m-i}$ \hs for \hs $0 \leq i \leq m$ \hs such that $\pi$
(respectively its nilradical $N$) consists of all \hs $x \in L$ \hs such that
\hs $xF_i \subset F_i$ \hs (resp. \hs $xF_i \subset F_{i-1}$) for \hs $1 \leq i
\leq m$.  \hs Put \hs $p_i = \dim (F_i/F_{i-1})$ \hs then it follows that \hs
$p_i = p_{m+1-i}$ \hs for \hs $1 \leq i \leq m$. \hs Following [B2, p.200] we
can find a Witt basis of $V$:
\begin{eqnarray*}
e_1,\ldots, e_r, \ e_{-r},\ldots, e_{-1}
\end{eqnarray*}
compatible with the given flag and such that \hs $\varphi(e_i, e_{-j}) =
\delta_{ij}$.\\
We now identify each \hs $x \in L$ \hs with its matrix with respect to this
basis, i.e. \hs $x = \left(\begin{array}{cc}A &B\\C &D\end{array}\right)$ \hs
where \hs $A, B, C, D$ \hs are \hs $r \times r$ \hs matrices such that \hs $B =
\widehat{B}$,\hs $C = \widehat{C}$, \hs $D = - \widehat{A}$, \hs where the
transformation $\widehat{\ }$ is the transpose relative to the second diagonal.
If \hs $x \in N$\hs then $x$ is of the form as shown in figure 2.
\begin{center}
\includegraphics{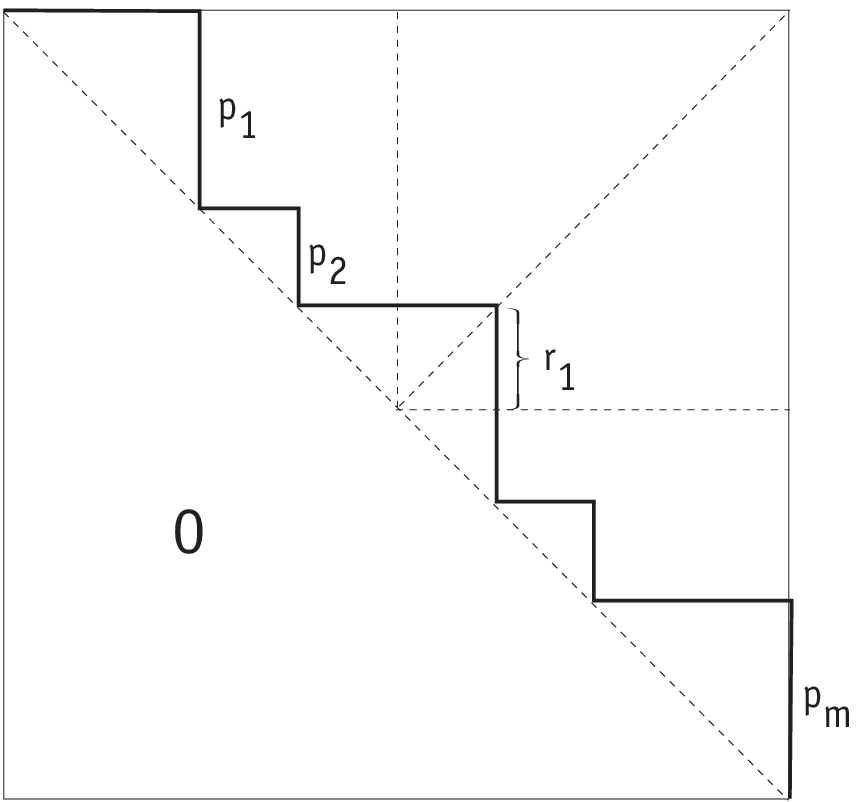}
\end{center}
\begin{center}
{\bf Figure 2}
\end{center}

If \hs $m = 2\ell + 1$ \hs then we put \hs $r_1 = \frac{1}{2} p_{\ell+1}$ \hs
($p_{\ell+1}$ is even since \hs $\ts\sum\limits_{i=1}^m p_i = n = 2r$ \hs and
\hs $p_i = p_{m+1-i}$). If \hs $m = 2\ell$ \hs then we put \hs $r_1 = 0$. \hs
$\pi$ is determined by the sequence \hs $(p_1,\ldots, p_\ell;r_1)$.  \hs Note
that \hs $r = \ts\sum\limits_{i=1}^\ell p_i + r_1$. \hs Next, let $P$ be the
subspace of $N$ of matrices of the form \hs $\left(\begin{array}{cc}
0 &B\\0 &0\end{array}\right)$ \hs where $B$ is an \hs $r \times r$ \hs matrix
such that \hs $B = \widehat{B}$ \hs and with zero \hs $r_1 \times r_1$ \hs
submatrix in the bottom left corner.  Clearly,
$$
X_{\varepsilon_i + \varepsilon_j} = E_{i,-j} + E_{j,-i}, \hs 1 \leq i \leq r-
r_1; \hs i \leq j \leq r
$$
form a basis of $P$ which is an abelian ideal of $N$ and \hs $\dim P =
\frac{1}{2} [(r^2 - r_1^2) + (r-r_1)]$. \hs We enlarge this basis to a basis of
$N$ by adjoining some vectors of the type
\begin{eqnarray*}
X_{\varepsilon_i - \varepsilon_j} = E_{ij} - E_{-j,-i}, \hs i < j
\end{eqnarray*}
From figure 2 we see that
\begin{eqnarray*}
\dim N &=& \ts\frac{1}{2} \left(r^2 - \ts\sum\limits_{i=1}^\ell p_i^2 -
r_1^2\right) + \dim P\\
&=& (r^2 - r_1^2) - \ts\frac{1}{2} \ts\sum\limits_{i=1}^\ell p_i^2 +
\ts\frac{1}{2} (r-r_1)
\end{eqnarray*}
Next, let \hs $f \in N^\ast$ \hs be defined by \hs $f(X_{2\varepsilon_i}) = 1$
\hs for \hs $1 \leq i \leq r-r_1$ \hs and zero on all other basis vectors of
$N$.  We want to show that \hs $P^f = P$.  \hs For this purpose we take \hs $x
\in P^f$ \hs which we can write as
\begin{eqnarray*}
x = \ts\sum\limits_{i < j} \lambda_{ij} X_{\varepsilon_i - \varepsilon_j} + y
\end{eqnarray*}
where \hs $X_{\varepsilon_i - \varepsilon_j} \in N$, \hs $\lambda_{ij} \in k$
\hs and \hs $y \in P$. \hs Fix any \hs $\lambda_{st}$, \hs $s < t$ \hs with \hs
$X_{\varepsilon_s - \varepsilon_t} \in N$. \hs This implies that \hs $s \leq r -
r_1$, \hs $t \leq r$. \hs Hence \hs $X_{\varepsilon_s + \varepsilon_t} \in P$.
\hs Therefore,
\begin{eqnarray*}
0 &=& f([x, X_{\varepsilon_s + \varepsilon_t}]) = \ts\sum\limits_{i < j}
\lambda_{ij} f([X_{\varepsilon_i - \varepsilon_j}, X_{\varepsilon_s +
\varepsilon_t}]) + f([y, X_{\varepsilon_s + \varepsilon_t}])\\
&=& \ts\sum\limits_{i < j} \lambda_{ij} f(\delta_{js} X_{\varepsilon_i +
\varepsilon_t} + \delta_{jt} X_{\varepsilon_i + \varepsilon_s})\\
&=& \ts\sum\limits_{i < s} \lambda_{is} f(X_{\varepsilon_i + \varepsilon_t}) +
\ts\sum\limits_{i < t} \lambda_{it} f(X_{\varepsilon_i + \varepsilon_s})\\
&=& 0 + \lambda_{st}
\end{eqnarray*}
($f(X_{\varepsilon_i + \varepsilon_t}) = 0$ \hs since \hs $i < s < t$).\\
It follows that \hs $x = y \in P$. \hs So, \hs $P^f \subset P$. \hs
Consequently, \hs $P^f = P$ \hs as the other inclusion is obvious.  By Lemma 1.2
we may conclude that $P$ is a CP of $N$ and \hs $f \in N_{\reg}^\ast$. \hs
Finally,
\begin{eqnarray*}
i(N) &=& 2\dim P - \dim N\\
&=& (r^2 - r_1^2) + (r-r_1) - (r^2 - r_1^2) + \ts\frac{1}{2}
\ts\sum\limits_{i=1}^\ell p_i^2 - \ts\frac{1}{2}(r-r_1)\\
&=& \ts\frac{1}{2} \left(\ts\sum\limits_{i=1}^\ell p_i^2 + (r-r_1)\right)\\
&=& \ts\frac{1}{2} \ts\sum\limits_{i=1}^\ell p_i(p_i + 1)
\end{eqnarray*}

\

\begin{center}
{\bf 7. CP-PRESERVING EXTENSIONS}
\end{center}

{\bf PROPOSITION 7.1} \hs Let $M$ be a finite dimensional Lie algebra over $k$
and let \hs $d \in \Der M$ \hs be a derivation such that \hs $d(Z(M)) \neq 0$.
\hs Consider the extension \hs $L = M \oplus kd$ \hs in which \hs $[d,x] =
d(x)$, $x \in M$.\\
Then we have
\begin{itemize}
\item[(i)] $i(M) = i(L) + 1$.
\item[(ii)] $L$ has a CP if and only if $M$ has a CP.
\item[(ii)] If $L$ is square integrable, then so is $M$.
\end{itemize}

{\bf REMARK 7.2} \hs Example (3) of 1.7 shows that the condition on $d$ cannot
be removed.

\

 {\it Proof.} \hs Take \hs $u \in Z(M)$ \hs such that \hs $d(u)
\neq 0$. \hs Clearly \hs $M = C(u)$.  \hs Now the assertions
follow directly from Proposition
1.8.

\

{\bf PROPOSITION 7.3} \hs Let $M$ be a finite dimensional Lie algebra over $k$
and fix $z$, a nonzero central element of $M$.  Let $S$ be a $2r$-dimensional
vector space, provided with a nondegenerate alternating bilinear form \hs
$\varphi : S \times S \rightarrow k$. \hs Consider the Lie algebra \hs $L = M
\oplus S$ \hs containing $M$ as an ideal and in which \hs $[x,s] = 0$ \hs and
\hs $[s,t] = \varphi(s,t)z$ \hs for \hs $x \in M$; \hs $s,t\in S$. \hs Then we
have
\begin{itemize}
\item[(i)] $H = S \oplus kz$ \hs is a Heisenberg Lie algebra
\item[(ii)] $i(L) = i(M)$ \hs and \hs $Z(L) = Z(M)$
\item[(iii)] $M$ is square integrable if and only if $L$ is square integrable
\item[(iv)] If $M$ allows a CP (resp. a CP-ideal) then the same holds for $L$.
\end{itemize}

{\it Proof.} \hs (i) It is easy to verify that $L$ is a Lie algebra.  There
exists a \hs $f \in L_{\reg}^\ast$ \hs such that \hs $f|_M \in M_{\reg}^\ast$
\hs and \hs $f(z) \neq 0$. \hs We may assume that \hs $f(z) = 1$ \hs (by
replacing $f$ by \hs $\ds\frac{1}{f(z)}f$). \hs Then for all $s,t \in S$
\begin{eqnarray*}
B_f(s,t) = f([s,t]) = \varphi(s,t)
\end{eqnarray*}
From the assumption on $\varphi$, \hs $S \ \cap \ S^\perp = 0$
and we can find a basis  $s_1,\ldots, s_r; t_1,\ldots,t_r$ \hs of
$S$ such that for all $i, j$:
\begin{eqnarray*}
\varphi(s_i,s_j) = 0 = \varphi(t_i, t_j)\hs \mbox{and} \hs \varphi(s_i,t_j) =
\delta_{ij}
\end{eqnarray*}
This implies \hs $[s_i, s_j] = 0 = [t_i, t_j]$ \hs and \hs $[s_i, t_j] =
\delta_{ij}z$ \hs for all $i, j$.  Consequently, $H$ is a Heisenberg Lie
algebra.\\
(ii) First, we notice that \hs $M = S^\perp$. \hs Indeed, \hs $M \subset
S^\perp$ \hs since \hs $f([M, S]) = 0$.  \hs For the other inclusion, take \hs
$x \in S^\perp$, \hs which we decompose as \hs $x = m+s$ \hs with \hs $m \in M$
\hs and \hs $s \in S$. \hs Then, \hs $s = x - m \in S \ \cap \ S^\perp = \{0\}$.
\hs Hence, \hs $x = m \in M$. \hs As \hs $M = S^\perp$ \hs we deduce from
[D,1.12.4] that
\begin{eqnarray*}
M(f|_M) = M \cap M^\perp = S \cap S^\perp + L^\perp = L(f)
\end{eqnarray*}
Taking dimensions yields \hs $i(M) = i(L)$. \hs Clearly, the elements of $Z(M)$
commute with those of $M$ and $S$.  Hence, \hs $Z(M) \subset Z(L)$. \hs
Conversely, take \hs $x \in Z(L)$ \hs which we can decompose as \hs $x = m+s$
\hs with \hs $m \in M$ \hs and \hs $s \in S$. \hs For all \hs $s' \in S$:
\begin{eqnarray*}
[s,s'] = [x-m, s'] = [x,s'] - [m, s'] = 0
\end{eqnarray*}
and hence also \hs $\varphi (s,s') = f([s,s']) = 0$ \hs which implies that \hs
$s=0$ \hs and so \hs $x = m \in M \cap Z(L) \subset Z(M)$.\\
(iii) This follows at once from (ii).\\
(iv) Suppose $P_1$ is a CP of $M$.  Put \hs $P_2 = ks_1 + \ldots + ks_r$ \hs and
\hs $P = P_1 \oplus P_2$. \hs Then $P$ is a CP of $L$ since $P$ is commutative
and
\begin{eqnarray*}
\dim P = \dim P_1 + \dim P_2 = \ts\frac{1}{2} (\dim M + i(M)) + \ts\frac{1}{2}
\dim S = \ts\frac{1}{2} (\dim L + i(L)).
\end{eqnarray*}
Finally, if $P_1$ is an ideal of $M$ then $P$ is an ideal of $L$ since
\begin{eqnarray*}
[M,P] = [M,P_1] + [M, P_2] = [M, P_1] \subset P_1 \subset P
\end{eqnarray*}
and
\begin{eqnarray*}
[t_j,P] &=& [t_j, P_1] + [t_j, P_2] = [t_j, P_2]\\
&=& \ts\sum\limits_i k[t_j, s_i] = kz \subset Z(M) \subset P_1 \subset P.
\end{eqnarray*}

{\bf PROPOSITION 7.4} \hs Let $A$ be an $n$-dimensional
commutative (associative) Fro\-be\-ni\-us algebra over $k$ and $M$
an $m$-dimensional Lie algebra
over $k$.\\
Consider the Lie algebra \hs $L = A \otimes_k M$ for which \hs $[a \otimes x, a'
\otimes y] = aa' \otimes [x,y]$, \hs $a,a' \in A$ \hs and \hs $x,y \in M$.  \hs
Then we have
\begin{itemize}
\item[(i)] $M$ is square integrable if and only if $L$ is square integrable.
\item[(ii)] $M$ is Frobenius if and only if $L$ is Frobenius.
\item[(iii)] If $M$ allows a CP (resp. a CP-ideal) then the same holds for $L$.
\end{itemize}

{\it Proof.} \hs (i) From [F, p.241-243] we know that \hs $i(L) = n.i(M)$. \hs
On the other hand, \hs $Z(L) = A \otimes_k Z(M)$ \hs and so \hs $\dim Z(L) =
n.\dim Z(M)$. \hs Therefore, \\
$i(L) = \dim Z(L)$ \hs if and only if \hs $i(M) = \dim Z(M)$.\\
(ii) This follows from (i) and its proof.\\
(iii) Let $P$ be a CP (resp. a CP-ideal) of $M$.  Then \hs $Q = A \otimes_k P$
\hs is a commutative Lie subalgebra (resp. ideal) of $L$ and
\begin{eqnarray*}
\dim Q = n.\dim P = n.\ts\frac{1}{2} (\dim M + i(M))\\
= \ts\frac{1}{2}(n.\dim M + n.i(M)) = \frac{1}{2} (\dim L + i(L))
\end{eqnarray*}

\

\begin{center}
{\bf ACKNOWLEDGMENTS}
\end{center}
The authors thank Jacques Alev for his interest and his valuable comments on the
subject.  They also express their gratitude to Dmitri Panyushev for providing
some useful preprints of his work.  The first author is grateful for the
generous hospitality he received during his visits at the University of Limburg
and the University of Bochum, which greatly contributed to the completion of
this paper.\\
Finally, the authors thank the referee for pointing out Remark
6.3(b) and for giving helpful suggestions in connection with the
presentation of the paper.

\bibliographystyle{amsalpha}

\end{document}